\documentclass[reqno,12pt]{amsart}
\usepackage{a4wide}
\usepackage{amssymb}
\usepackage{amsmath}
\usepackage[pctex32]{graphics}
\usepackage{graphicx}

\newtheorem{theorem}{Theorem}

\newtheorem{lemma}[theorem]{Lemma}
\newtheorem{proposition}[theorem]{Proposition}

\theoremstyle{definition}

\newtheorem{remark}[theorem]{Remark}

\newcommand{\ds}{\displaystyle}
\newcommand{\dis}{\displaystyle}
\newcommand{\be}{\begin{equation}}
\newcommand{\bd}{\begin{displaymath}}
\newcommand{\ed}{\end{displaymath}}
\newcommand{\ba}{\begin{eqnarray}}
\newcommand{\ea}{\end{eqnarray}}

\newcommand{\p}{\partial}

\def\Realr{\mathbb{R}}
\def\rr{\mathbb R}
\def\ee{\mathbb E}
\def\eps{\varepsilon}

\def\t{\tau}

\newcommand{\RN}{\mathbb{R}^{3}}
\renewcommand{\P}{{\mathcal P}}               
\newcommand{\Ptwo}{\P_2(\RN)}             

\def\proof{\par{\it Proof.-} \ignorespaces}
\def\endproof{{\ \vbox{\hrule\hbox{%
   \vrule height1.3ex\hskip0.8ex\vrule}\hrule
  }}\par}

\numberwithin{equation}{section}

\begin{document}

\title{Tanaka Theorem for Inelastic Maxwell Models}

\author{Fran\c cois Bolley}
\address{Institut de Math\'ematiques, LSP (UMR C5583), Universit\'e Paul Sabatier, Route de Nar\nobreak bonne, F-31062 Toulouse cedex 9}
\email{bolley@cict.fr}

\author{Jos\'e A. Carrillo}

\address{ICREA (Instituci\'o Catalana de Recerca i
Estudis Avan\c cats) and Departament de Mate\-m\`a\-tiques,
Universitat Aut\`onoma de Barcelona, E-08193 Bellaterra, Spain}
\email{carrillo@mat.uab.es}


\maketitle


\begin{abstract}
We show that the Euclidean Wasserstein distance is contractive for
inelastic homogeneous Boltzmann kinetic equations in the
Maxwellian approximation and its associated Kac-like caricature.
This property is as a generalization of the Tanaka theorem to
inelastic interactions. Even in the elastic classical Boltzmann
equation, we give a simpler proof of the Tanaka theorem than
the ones in \cite{Tanaka,Vil03}.  Consequences are drawn
on the asymptotic behavior of solutions in terms only of the Euclidean
Wasserstein distance.
\end{abstract}

\section{Introduction}
This work is devoted to contraction and asymptotic properties of
the homogeneous Boltzmann-type equations for inelastic
interactions in the Maxwe\-llian approximation introduced
in~\cite{Bobylev-Carrillo-Gamba} and further analyzed in
\cite{Carrillo-Cercignani-Gamba,Bobylev-Cercignani,Bobylev-Cercignani2,Bobylev-Cercignani-Toscani,BisiCT,Bobylev-Gamba,BisiCT2,Bobylev-Cercignani-Gamba}.
We are basically concerned with the Boltzmann equation
\begin{equation}
\frac{\p f}{\p t}=B \sqrt{\theta(f(t))} \, Q(f,f) \label{BE1}
 \end{equation}
considered in~\cite{Bobylev-Carrillo-Gamba} and its variants.
Here, $f(t,v)$ is the density for the velocity $v\in\RN$
distribution of the molecules at time $t$, and $Q(f,f)$ is the
inelastic Boltzmann collision operator defined by
\begin{equation}
(\varphi,Q(f,f))= \frac{1}{4 \pi} \int_{\Realr^3} \int_{\Realr^3}
\int_{S^2} f(v) f(w) \Big[ \varphi (v')- \varphi (v) \Big]
d\sigma\,dv\,dw \label{Qweak}
 \end{equation}
for any test function ~$\varphi$, where
\[
 \dis v'= \frac{1}{2} (v+w)+ \frac{1-e}{4}(v-w)+\frac{1+e}{4} |v-w|\sigma
\]
is the postcollisional velocity, $\sigma\in S^2$, $v, w \in \rr^3$  and $0 < e  \leq
1$ is the constant restitution coefficient. Equation \eqref{BE1}
preserves mass and momentum, but makes the kinetic energy (or temperature)
$$
\theta(f(t))=\frac13 \int_{\RN} \Big\vert v - \int_{\rr^3} v \, f(t, v) \, dv \Big|^2 f(t,v)\,dv
$$
decrease towards $0$. In particular, solutions to \eqref{BE1} tend
to the Dirac mass at the mean velocity of the particles
\cite{Bobylev-Carrillo-Gamba}. We refer to
\cite{Bobylev-Carrillo-Gamba,Bobylev-Cercignani2,Vil06} for the
discussion about the relation of this model to the inelastic
hard-sphere Boltzmann equation and different ways of writing the
operator. Let us just point out that the factor $B \,
\sqrt{\theta(f(t))}$ in front of the operator in \eqref{BE1} is
chosen for having the same temperature decay law as its
hard-sphere counterpart \cite{Bobylev-Carrillo-Gamba} known as the
Haff's law.

The convergence towards the monokinetic distribution has been made
more precise in
\cite{Bobylev-Cercignani2,Bobylev-Cercignani-Toscani,BisiCT2} by
means of homogeneous cooling states. They are self-similar
solutions of the homogeneous Boltzmann equation \eqref{BE1}
describing the long-time asymptotics and presenting power-like
tail behavior whose relevance was previously discussed in the
physics literature \cite{EB1,EB2}.

To avoid the collapse of the solution to the Dirac mass, the
authors in \cite{Carrillo-Cercignani-Gamba} suggested the
introduction of a stochastic thermostat which, at the kinetic
level, is modelled by a linear diffusion term in velocity. In this
framework, the density $f$ in the velocity space obeys
\begin{equation}
\frac{\p f}{\p t}=B \sqrt{\theta(f(t))}\, Q(f,f) + A\, \theta^p
(f(t)) \, \triangle_v f \qquad \,\,\, {\rm with} \quad \,\,\, 0 \leq
p < \frac{3}{2} \cdot \label{BE2}
\end{equation}
Existence and uniqueness for given mean velocity of a steady state
to \eqref{BE2} have been shown in
\cite{Cercignani-Illner-Stoica,Bobylev-Cercignani,BisiCT}. The
convergence of solutions towards this steady state in all Sobolev
norms has also been investigated and quantified by means of
Fourier-based distances between probability measures
\cite{BisiCT}.

Fourier techniques are a good toolbox and have been extremely
fruitful for studying Maxwellian models in kinetic theory since
Bobylev observed \cite{BobylevBKW,Bobylev} that such equations
have closed forms in Fourier variables. Fourier distances are not
only suitable technical tools to study the long-time asymptotics
of models \eqref{BE1} and \eqref{BE2}, but also they represent the
first Liapunov functionals known for inelastic Boltzmann-type
equations \cite{BisiCT,BisiCT2}. In the case of the classical
elastic Boltzmann equation for Maxwellian molecules, there is
another known Liapunov functional, namely, the Tanaka functional
\cite{Tanaka}, apart from the $H$-functional for which no
counterpart is known in inelastic models.

The Tanaka functional is the Euclidean (or quadratic) Wasserstein distance
between measures in the modern jargon of optimal mass transport
theory. It is defined on the set $\Ptwo$ of Borel probability measures
on $\rr^3$ with finite second moment or kinetic energy  as
\[
W_2(f, g) = \inf_{\pi} \left\{ \iint_{\rr^3 \times \rr^3} \vert v
-w \vert^2 \, d\pi(v, w) \right\}^{1/2}= \inf_{(V,W)} \left\{
\ee\left[ \vert V-W \vert^2 \right] \right\}^{1/2}
\]
where $\pi$ runs over the set of joint probability measures on
$\rr^3 \times \rr^3$ with marginals $f$ and $g$ and $(V,W)$ are
all possible couples of random variables with $f$ and $g$ as
respective laws. This functional was proven by Tanaka
\cite{Tanaka} to be non-increasing for the flow of the homogeneous
Boltzmann equation in the Maxwellian case. In fact, the Tanaka
functional and Fourier-based distances are related to each other
\cite{Gabetta-Toscani-Wennberg,CGT,Toscani-Villani}, and were used
to study the trend to equilibrium for Maxwellian gases. On the
other hand, related simplified granular models
\cite{Carrillo-McCann-Villani} have been shown to be strict
contractions for the Wasserstein distance $W_2$.

With this situation, a natural question arose as an open problem
in \cite[Remark 3.3]{BisiCT2} and \cite[Section 2.8]{Vil06}: is
the Euclidean Wasserstein distance a contraction for the flow of
inelastic Maxwell models? The main results of this work answer
this question affirmatively. Moreover, we shall not need to introduce
Bobylev's Fourier representation of the inelastic Maxwell models
working only in the phy\-sical space.

We shall show in the next section the key idea behind the proof of all
results concerning contractions in $W_2$ distance for inelastic Maxwell
models, namely, the gain part $Q^+(f,f)$ of the collision operator
verifies
\[
W_2(Q^+(f, f) , Q^+(g, g) ) \leq \sqrt{\frac{3+ e^2}{4}} \, W_2(f,
g)
\]
for any $f,g$ in $\Ptwo$ with equal mean velocity and any restitution
coefficient $0<e\leq~1$. Based on this property, we shall
derive contraction and asymptotic properties both for \eqref{BE1}
and \eqref{BE2} in Subsections \ref{sec-nodiff} and
\ref{sec-diff}. On one hand, we shall prove that the flow for the diffusive
equation \eqref{BE2} is a strict contraction for $W_2$, while for
the scaled equation associated to \eqref{BE1} we shall show
that solutions converge in $W_2$ to a corresponding homogeneous cooling
state, without rate but only assuming that initial data
have bounded second moment. This improves the Ernst-Brito
conjecture
\cite{EB1,EB2,Bobylev-Cercignani2,Bobylev-Cercignani-Toscani,BisiCT2}
since it shows that the basin of attraction of the homogenous
cooling state is larger -we avoid the typical assumption of
bounded moments of order $2+\delta$- if we do not ask for a rate.

Moreover, a generalization for non constant cross sections
including Tanaka's theorem as a particular case will be proven in
Section \ref{sec-gene-b}. Finally, we shall also show this generic
property for the inelastic Kac model introduced in \cite{PT03} as
a dissipative version of Kac's caricature of Maxwellian gases
\cite{Kac,McKean}.


\section{Contraction in $W_2$ of the gain operator}\label{key}

We start by summarizing the main properties of the Euclidean
Wasserstein distance $W_2$ that we shall make use of in the rest,
 refering to \cite{Bolley,Vil03} for the proofs.

\begin{proposition}\label{w2properties}
The space $(\Ptwo,W_2)$ is a complete metric space. Moreover, the
fol\-lowing properties of the distance $W_2$ hold:
\begin{enumerate}
\item[i)] {\bf Convergence of measures:} Given $\{f_n\}_{n\ge 1}$
and $f$ in $\Ptwo$, the following three assertions are equivalent:
\begin{itemize}
\item[a)] $W_2(f_n, f)$ tends to $0$ as $n$ goes to infinity.

\item[b)] $f_n$ tends to $f$ weakly-* as measures as $n$ goes to infinity and
$$
\sup_{n\ge 1} \int_{\vert v \vert > R} \vert v \vert^2 \, f_n(v)
\, dv \to 0 \, \mbox{ as } \, R \to +\infty.
$$

\item[c)] $f_n$ tends to $f$ weakly-* as measures and
$$
 \int_{\rr^3} \vert v \vert^2 \, f_n(v) \, dv \to
\int_{\rr^3} \vert v \vert^2 \, f(v) \, dv  \, \mbox{ as } \, \mbox{n} \to + \infty.
$$
\end{itemize}

\item[iii)] {\bf Relation to Temperature:} If $f$ belongs to $\Ptwo$ and
$\delta_a$ is the Dirac mass at $a$ in $\RN$, then
$$
W_2^2(f,\delta_a)=\int_{\RN} |v-a|^2 df(v).
$$

\item[iii)] {\bf Scaling:} Given $f$ in $\Ptwo$ and $\theta>0$, let
us define
$$
{\mathcal S}_\theta [f]=\theta^{3/2} f(\theta^{1/2}v)
$$
for absolutely continuous measures with respect to Lebesgue measure or its
corres\-ponding definition by duality for general measures; then for
any $f$ and $g$ in $\Ptwo$, we have
$$
W_2({\mathcal S}_\theta [f],{\mathcal S}_\theta [g]) = \theta^{-1/2}\,
W_2(f,g) .
$$

\item[iv)] {\bf Convexity:} Given $f_1$, $f_2$, $g_1$ and $g_2$ in
$\Ptwo$ and $\alpha$ in $[0,1]$, then
$$
W_2^2(\alpha f_1 + (1-\alpha) f_2,\alpha g_1 + (1-\alpha) g_2)
\leq \alpha W_2^2(f_1,g_1) + (1-\alpha) W_2^2(f_2,g_2).
$$
As a simple consequence, given $f,g$ and $h$ in $\Ptwo$, then
$$
W_2(h * f,h * g) \leq W_2(f,g)
$$
where $*$ stands for the convolution in $\rr^3$.
\end{enumerate}
\end{proposition}

\medskip

\medskip

\noindent
Here the convolution of the two measures $h$ and $f$ is defined by duality by
$$
(\varphi, h * f) = \iint_{\rr^3 \times \rr^3} \varphi(x+y) \, dh(x) \, df(y)
$$
for any test function $\varphi$ on $\rr^3$. If $f$ is a Borel probability measure on $\rr^3$ we
shall let
\[
<\!\!f\!\!> = \int_{\rr^3} v \, df(v) =  \int_{\rr^3} v \,
f(v)\,dv
\]
denote its mean velocity. We shall use the same notation for
densities and measures expecting that the reader will not get
confused.

\medskip

Let us write the collision operator $Q$ given in \eqref{BE2} as
\begin{equation}\label{decompoQ}
Q(f,f) = Q^+(f, f) - f
\end{equation}
where $Q^+(f, f)$ is defined by
\begin{equation}\label{Q+weak}
(\varphi,Q^+(f,f))= \frac{1}{4 \pi} \int_{\Realr^3}
\int_{\Realr^3} \int_{S^2} f(v) f(w) \, \varphi (v') \, d\sigma\,dv\,dw
\end{equation}
for any test function $\varphi$, where we recall that
\[
 \dis v'= \frac{1}{2} (v+w)+ \frac{1-e}{4}(v-w)+\frac{1+e}{4} |v-w|\sigma.
\]

In this section we derive a contraction property in $W_2$ distance
of the gain operator $Q^+$. For that purpose, let us note that the
previous definition of the gain operator can be regarded as
follows: given a probability measure $f$ on $\rr^3$, the
probability measure $Q^+(f, f)$ is defined by
\[
 (\varphi,Q^+(f,f))= \int_{\Realr^3} \int_{\Realr^3} f(v) \, f(w) \, (\varphi, \Pi_{v, w}) \, dv \, dw
\]
where $\Pi_{v, w}$ is the uniform probability
distribution on the sphere $S_{v, w}$ with center
\[
 c_{v, w} = \frac{1}{2} (v+w)+ \frac{1-e}{4}(v-w)
\]
and radius
\[
 r_{v, w} = \frac{1+e}{4} |v-w|.
\]

\begin{figure}[ht]
\centering
\includegraphics[scale=0.7]{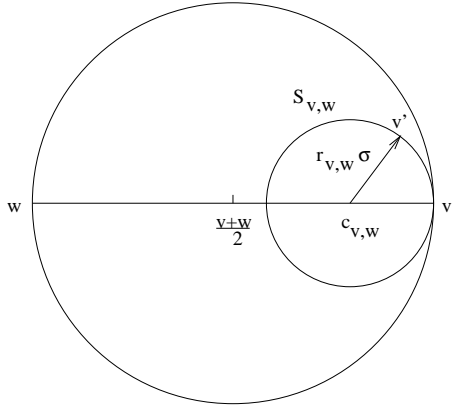}
\caption{Geometry of inelastic collisions} \label{spherepicture}
\end{figure}

In probabilistic terms, the gain operator is defined as an
expectation:
$$
Q^+(f, f)=\ee \left[ \Pi_{V, W}\right]
$$
where $V$ and $W$ are independent
random variables  with law $f$.

Then the convexity of $W_2^2$ in Proposition \ref{w2properties} implies
\begin{align}
W_2^2(Q^+(f,f), Q^+(g,g)) &= W_2^2(\ee \left[ \Pi_{V, W}\right],
\ee \left[ \Pi_{X, Y}\right]) \nonumber \\
&\leq \ee \left[ W_2^2(\Pi_{V, W}, \Pi_{X, Y})\right]
\label{conv-ee}
\end{align}
where $X$ and $Y$ are independent random variables with law $g$.
This observation leads us to consider the $W_2$ distance between
uniform distributions on spheres. To this aim, we have the
following general lemma:

\begin{lemma}\label{w2unif}
The squared Wasserstein distance $W_2^2$ between the uniform
distributions on the sphere with center $O$ and radius $r$
and the sphere with center $O'$ and radius $r'$ in $\rr^3$
is bounded by $\vert O'-O \vert^2 + (r' - r)^2$.
\end{lemma}

\proof We define a map $T:\RN\longrightarrow\RN$ transporting the
sphere of center $O$ and radius $r>0$ onto the sphere with center
$O'$ and radius $r'\ge r$ in the following way:
\begin{itemize}
\item If $r=r'$, then we just let $T$ be the translation map with vector
  $O'-O$, i.e., $T(v)=v+O'-O$.

\item If $O=O'$, then we just let $T$  be the dilation with factor
  $\frac{r'}{r}$ centered at $O$, i.e., $T(v)=\frac{r'}{r}v$.

\item If $r \neq r'$, then we consider the only point $\Omega\in
\RN$ verifying that
$$
\frac{1}{r} (O - \Omega) = \frac{1}{r'} (O' - \Omega),
$$
that is,
$$
\Omega = O + \frac{r}{r'-r} (O' - O).
$$
Then we let $T$ be the dilation with factor $\frac{r'}{r}$
centered at $\Omega$, that is, we let $T(v)=\Omega + \frac{r'}{r}
(v-\Omega)$. Such a construction of the point $\Omega$ and the map
$T$ is sketched in Figure \ref{lemsphere} in the case of non
interior spheres.

\begin{figure}[ht]
\centering
\includegraphics[scale=0.7]{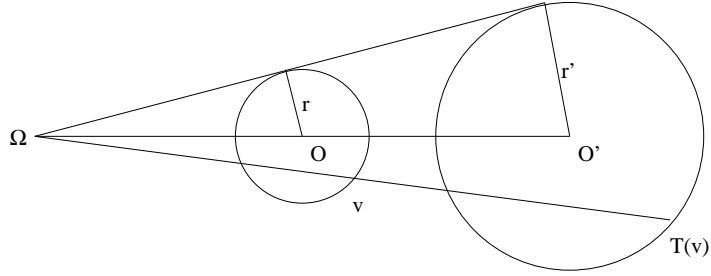}
\caption{Sketch of the computation of the Euclidean cost of
transporting spheres to spheres. Transport lines are just rays
from the point $\Omega$.} \label{lemsphere}
\end{figure}

\end{itemize}

Let ${\mathcal U}_{O,r}$ and ${\mathcal U}_{O',r'}$ denote the
uniform distributions on the corresponding spheres. Then the
transport plan $\pi$ given by
$$
\iint_{\RN\times\RN} \eta(v,w)\, d\pi(v,w) = \int_{\RN}
\eta(v,T(v))\, d{\mathcal U}_{O,r}(v)
$$
for all test functions $\eta(v,w)$ has ${\mathcal U}_{O,r}$ and ${\mathcal U}_{O',r'}$
as marginals by construction of $T$.
Using this transference plan in the definition of the Euclidean
Wasserstein distance, we finally conclude
$$
W_2^2({\mathcal U}_{O,r},{\mathcal U}_{O',r'}) \leq \int_{\RN} |v-T(v)|^2 \,
d{\mathcal U}_{O,r}(v)\! = \!\left(\frac{r'-r}{r}\right)^2 \int_{\RN}
|v-\Omega|^2 \, d{\mathcal U}_{O,r}(v)
$$
that can be computed explicitly, giving
$$
W_2^2({\mathcal U}_{O,r},{\mathcal U}_{O',r'}) \leq |O'-O|^2 + (r'-r)^2
$$
and finishing the proof.
\endproof

\medskip

This lemma, using the notation $a=v-x$ and $b=w-y$, for fixed
values $v, w, x, y$ in $\RN$, implies that
\begin{eqnarray*}
W_2^2(\Pi_{v, w}, \Pi_{x, y}) &\leq&
\vert c_{v, w} - c_{x, y} \vert^2 + \vert r_{v, w} - r_{x, y} \vert^2 \\
&\leq&
\Big\vert \frac{3-e}{4} \, a + \frac{1+e}{4} \, b \Big\vert^2 +
\Big( \frac{1+e}{4} \Big)^2 \, \vert a - b \vert^2 \\
&=& \frac{5 - 2 \, e + e^2}{8} \, \vert a \vert^2 + \frac{(1 +
e)^2}{8} \, \vert b \vert^2 + \frac{1-e^2}{4} \, a\cdot b;
\end{eqnarray*}
here $a \cdot b$ denotes the scalar product between $a$ and $b$ in
$\rr^3$ and the bound in
\begin{align*}
\vert r_{v, w} - r_{x, y} \vert^2  &= \Big( \frac{1+e}{4} \Big)^2
\, \big\vert \vert v-w \vert - \vert x-y \vert \big\vert^2 \\
& \leq \Big( \frac{1+e}{4} \Big)^2 \, \big\vert (v-w) - (x-y)
\big\vert^2 = \Big( \frac{1+e}{4} \Big)^2 \, \vert a - b
\vert^2
\end{align*}
follows from the Cauchy-Schwarz inequality
\begin{equation}\label{csineq}
(v-w) \cdot (x-y) \leq \vert v-w \vert \, \vert x-y \vert.
\end{equation}
Therefore, by \eqref{conv-ee},
\begin{align*}
W_2^2(Q^+(f,f), Q^+(g,g))
\leq  & \;
\frac{5 - 2 \, e + e^2}{8} \, \ee \left[ \vert V-X \vert^2 \right]
+\frac{(1 + e)^2}{8} \, \ee \left[
\vert W-Y \vert^2 \right] \\
& + \, \frac{1-e^2}{4} \, \ee \left[ (V-X) \cdot
(W-Y)\right].
\end{align*}

Let moreover  $(V, X)$ and $(W, Y)$ be two independent optimal couples in
the sense that
\[
W_2^2(f, g) = \ee \left[ \vert V-X \vert^2\right] =  \ee \left[
\vert W-Y \vert^2\right].
\]
Then
\[
\ee \left[ (V-X) \cdot (W-Y)\right] =  \ee \left[ (V-X) \right]\,  \cdot
\, \ee \left[ (W-Y)\right]  = \big\vert <\!f\!> - <\!g\!> \big\vert^2
\]
 by independence. Collecting all terms leads to the following key estimate and
contraction property:

\begin{proposition}\label{prop-qplus}
If $f$ and $g$ belong to $\Ptwo$, then
\[
W_2^2(Q^+(f, f) , Q^+(g, g) ) \leq \frac{3+ e^2}{4} \, W_2^2(f, g)
+ \frac{1-e^2}{4} \, \big\vert <\!f\!> - <\!g\!> \big\vert^2
\]
for any restitution coefficient $0<e\leq 1$. As a consequence,
given $f$ and $g$ in $\Ptwo$ with equal mean velocity, then
\[
W_2(Q^+(f, f) , Q^+(g, g) ) \leq \sqrt{\frac{3+ e^2}{4}} \, W_2(f,
g).
\]
\end{proposition}

\medskip

The case of equality is addressed in the following statement:

\begin{proposition}\label{eg-Q+}
Let $f$ and $g$ belong to $\Ptwo$ with equal mean velocity and
temperature, where $g$ is absolutely continuous with respect to
Lebesgue measure with positive density. If
\[
W_2(Q^+(f, f) , Q^+(g, g) ) = \sqrt{\frac{3+ e^2}{4}} \, W_2(f, g).
\]
for some restitution coefficient $0<e\leq 1$, then $f=g$.

\end{proposition}

\proof It is necessary that the equality holds at each step of the
arguments in Proposition~\ref{prop-qplus}. In particular,
\eqref{csineq} holds as an equality, that is,
$$
\frac{V - W}{\vert V - W \vert} = \frac{X-Y}{\vert X-Y \vert}
$$
almost surely in the above notation. Then, since $g$ is absolutely
continuous with respect to Lebesgue measure with positive density,
one can proceed as in \cite[Lemma 9.1]{Tanaka} to show that $f=g$.
We sketch the proof for the sake of the reader. Since $g$ is
absolutely continuous with respect to Lebesgue measure, there
exists \cite{Vil03} a Borel map $u: \rr^3 \to \rr^3$ such that $f$
be the image measure of $g$ by $u$, and in probabilistic terms $V
= u(X)$ and $W = u(Y)$ almost surely. Hence
\begin{equation}\label{cseq}
\frac{u(x) - u(y)}{\vert  u(x) - u(y) \vert} = \frac{x-y}{\vert
x-y \vert}
\end{equation}
almost everywhere for Lebesgue measure since $X$ and $Y$ are
independent and since their law $g$ has positive density. We leave
the reader to check \cite[Exercise~7.25]{Vil03} that this implies
the existence of constants $\omega_1$ and $\omega_2$ such that
$u(x)=\omega_1 + \omega_2 x$. First of all $\omega_2^2 =1$ since
$f$ and $g$ have same temperature. Then  identity \eqref{cseq}
forces $\omega_2 =1$, implying $\omega_1=0$ since $< \! f  \!> =
<\! g \! >$, and finally $f=g$.
\endproof


\section{Contractive Estimates for the Inelastic Maxwell Model}

In this section, we shall derive contractive estimates in the
Euclidean Wasserstein distance for solutions to the inelastic
Maxwell models both in the non-diffusive and the diffusive cases.

\subsection{The non-diffusive case}\label{sec-nodiff}

We are first concerned with solutions $f(t)$ to the Boltzmann equation
\eqref{BE1} with $0 < e < 1$. After time scaling defined by
\[
\t = \frac{B}{E} \int_0^t \sqrt{\theta(f(w))} \, dw
\]
with $\ds E = \frac{8}{1 - e^2} \raise2pt \hbox{,}$  as in
\cite{BisiCT2}, we get a function denoted again $f(\t)$ for
simplicity, solution to
\begin{equation}
\frac{\p f}{\p \t}= E \, Q(f,f) \label{BE1tau}.
\end{equation}

\begin{theorem}\label{thm-BE1tau}
If $f_1$ and $f_2$ are two solutions to \eqref{BE1tau} with
respective initial data $f_1^0$ and $f_2^0$ in $\Ptwo$, then
\begin{align}
W_2^2(f_1(\t), f_2(\t)) \leq {\rm e}^{-2 \tau}\, W_2^2(f_1^0,
f_2^0) + (1 - {\rm e}^{-2 \tau})\,\big\vert \!<\!f_1^0\!> -
<\!f_2^0\!> \!\big\vert^2 \label{cont-BE1tau}
\end{align}
for all $\t \geq 0$.
\end{theorem}

\proof Decomposition \eqref{decompoQ} of the
collision operator $Q$ as
$$
Q(f,f) = Q^+(f, f) - f
$$
 allows us to represent the solutions to \eqref{BE1tau} by Duhamel's formula
as
\[
f_i(\t) = {\rm e}^{-E \t} \, f_i^0 + E \int_0^{\t} {\rm
e}^{-E(\t-s)} \, Q^+(f_i(s), f_i(s)) \, ds, \qquad i=1,2.
\]
Then the convexity of the squared Wasserstein distance in Proposition
\ref{w2properties} and Proposition \ref{prop-qplus} imply
\begin{eqnarray*}
&& W_2^2(f_1(\t), f_2(\t)) \\
&& \leq {\rm e}^{-E \t} \, W_2^2(f_1^0, f_2^0) + E \int_0^{\t}
\!\! {\rm e}^{-E(\t-s)} \, W_2^2 \big(
Q^+(f_1(s), f_1(s)), Q^+(f_2(s), f_2(s)) \big) \, ds \\
&& \leq {\rm e}^{-E \t} \, W_2^2(f_1^0, f_2^0) + E \int_0^{\t}
\!\! {\rm e}^{-E(\t-s)} \left( \frac{3+e^2}{4} \, W_2^2 (f_1(s),
f_2(s)) + X \right) ds;
\end{eqnarray*}
here
$$
X=\frac{1-e^2}{4} \, \big\vert <\!\!f_1(s)\!\!> -
<\!\!f_2(s)\!\!> \big\vert^2
$$
does not depend on time since the mean velocity is preserved by
equation \eqref{BE1tau}. In other words, the function $y(\t) =
{\rm e}^{E \t} \, W_2^2(f_1(\t), f_2(\t))$ satisfies the
inequality
 \[
 y(\t) \leq y(0) + E \int_0^{\t} \left( \frac{3+e^2}{4} \, y(s) + X \, {\rm e}^{E s} \right) \, ds
 \]
and then
\[
 y(\t) \leq y(0) \, {\rm e}^{\gamma E \t} + \frac{X}{1 - \gamma} ({\rm e}^{E \t} - {\rm e}^{\gamma E \t})
\]
by Gronwall's lemma with $\gamma=(3+e^2)/4$. This concludes the
argument since $(1 - \gamma ) \, E = 2$. \endproof

\begin{remark}
\
\begin{enumerate}
\item Without further assumptions on the initial data $f_1^0$ and
$f_2^0$, this result is optimal in the following sense. If $f_2^0$
is chosen as the Dirac mass at the mean velocity of $f_1^0$, then
inequality \eqref{cont-BE1tau} is actually an equality for all
$\t$; indeed
\begin{align*}
W_2^2(f_1(\t), f_2(\t)) &= \int_{\rr^3} \vert v \, -
<\!\!f_1(\t)\!\!> \! \vert^2 \, f_1(\t, v) \, dv = 3\, \theta(f_1(\t)) \\
&= 3 \, {\rm e}^{-2 \tau} \, \theta(f_1^0) = {\rm e}^{-2 \, \tau}
\, W_2^2(f_1^0, f_2^0)
\end{align*}
since $\displaystyle \frac{d \theta}{d \t} = -\,2\, \theta$ by equation \eqref{BE1tau}.

\item In terms of the original time variable $t$ in \eqref{BE1}, if $f_1^0$ and
$f_2^0$ are two initial data with the same initial temperature
$\theta_0$, then the temperatures of the corresponding solutions $f^1$ and $f^2$
to \eqref{BE1} follow the law
\begin{equation}
\frac{d \theta}{dt} = -\, \frac{1-e^2}{4}\, B
\theta^{\frac{3}{2}} \label{eqtemp}
\end{equation}
and hence are both equal to
$$
\theta(t)= \left(\theta_0^{-1/2} + \frac{1-e^2}{8}B t\right)^{-2}.
$$
Then estimate \eqref{cont-BE1tau} reads as
$$
W_2^2(f_1(t), f_2(t)) \leq \frac{\theta(t)}{\theta_0} \,
W_2^2(f_1^0, f_2^0) + \left(1 - \frac{\theta(t)}{\theta_0} \right)
\big\vert\!<\!f_1^0\!> \!-\! <\!f_2^0\!> \big\vert^2
$$
for all $t \geq 0$.

\end{enumerate}
\end{remark}

\medskip

\noindent The convergence of the solutions to
\eqref{BE1} towards the Dirac measure at their mean velocity has been
made precise in \cite{Bobylev-Cercignani-Toscani,BisiCT2} by the
introduction of self-similar vari\nobreak ables and homogeneous cooling
states. There the authors prove that the rescaled solutions $g$
defined by
\begin{equation}\label{f-g}
g(\t, v) = \theta^{3/2} (f(\t)) \, f(\t, \theta^{1/2} (f(\t)) \, v)
\end{equation}
satisfy the strict contraction property
\[
d_{2 + \eps} (g_1(\t), g_2(\t)) \leq {\rm e}^{- C(\eps) \t} d_{2 +
\eps} (g_1^0, g_2^0), \qquad C(\eps) >0
 \]
for initial data $g_1^0$ and $g_2^0$ in $\Ptwo$ with equal mean
velocity and pressure tensor, where $\eps>0$ and $d_{2 + \eps}$ is a Fourier-based
distance between probability measures. Moreover, for $\eps =0$ one has $C(\eps) =0$ giving
a non-strict contraction in $d_2$ distance. In fact, by the scaling property in Proposition
\ref{w2properties}, \eqref{cont-BE1tau} reads as
\begin{equation}\label{contr1g}
 W_2(g_1(\t), g_2(\t)) \leq W_2 (g_1^0, g_2^0)
\end{equation}
in the scaled variables. This is consistent with the fact that the
distances $d_2$ and $W_2$ are ``of the same order''
\cite{Gabetta-Toscani-Wennberg,Toscani-Villani,BisiCT2} up to
moment bounds.

A measure $g(\t, v)$ defined by \eqref{f-g} from a
solution $f(v, \t)$ to \eqref{BE1tau} with initial zero mean
velocity has zero mean velocity and unit kinetic energy for all
$\t$, and is solution to
\begin{equation}\label{pde-g1}
\frac{\partial g}{\partial \t} + \nabla \cdot (g \, v) = E \, Q(g, g).
\end{equation}
Moreover it is proven in \cite{Bobylev-Cercignani-Toscani,BisiCT2}
that \eqref{pde-g1} has a unique stationary solution $g_{\infty}$
with zero mean velocity and unit kinetic energy; all measure
solutions $g(\t, v)$ to \eqref{pde-g1} with zero mean velocity,
unit kinetic energy and bounded moment of order $2 + \eps$ converge to
this stationary state $g_{\infty}$ as $\t$ goes to infinity
in the $d_2$ sense, that is, in the $W_2$ sense since $d_2$ and $W_2$
metrize the same topology on probability measures \cite{Toscani-Villani}
up to moment conditions. Moreover the convergence has exponential rate in the $d_2$
sense, and in the $W_2$ sense if the initial datum has finite fourth order moment.
In turn this ensures existence and uniqueness of homogeneous cooling states to \eqref{BE1}
for given mean velocity and kinetic energy, and algebraic convergence of the solutions
$f(t)$ towards them in the original variables.

\medskip

We conclude this section by proving this convergence
result using only the $W_2$ distance, and without assuming
that the initial data has more than two finite moments. This in turn shows that the
Euclidean Wasserstein distance $W_2$ between solutions of
\eqref{pde-g1} converges to zero as $t$ goes to infinity, improving over
\eqref{contr1g} that does not {\it a priori} yield any information
on the long-time behavior of the solutions $g$. As a drawback, this argument does not
provide any rate of convergence as does the Fourier-based argument in
\cite{BisiCT2}.

\begin{theorem}\label{norate}
Let $g^0_1$ and $g^0_2$ be two Borel probability measures on
$\rr^3$ with zero mean velocity and unit kinetic energy, and let
$g_1(\t)$ and $g_2(\t)$ be the solutions to \eqref{pde-g1} with
respective initial data $g_1^0$ and $g_2^0$. Then the map $\t \mapsto
W_2(g_1(\t), g_2(\t))$ is non-increasing and tends to $0$ as $\t$
goes to infinity.
\end{theorem}

\medskip

\proof It is based on the argument in
\cite{Vil03} to Tanaka's theorem. The first statement is a simple consequence of
\eqref{contr1g}. Then we turn to the second part of the
theorem which by triangular inequality for the $W_2$ distance
is enough to prove when $g_2^0$, and hence $g_2(\t)$, is the
unique stationary state $g_{\infty}$  to \eqref{pde-g1} with zero
mean velocity and unit kinetic energy.

\medskip
{\bf Step 1.-} Let us first assume that the fourth moment of the
initial datum is bounded, i.e.,
$$
\int_{\rr^3} \vert v \vert^4 \, g^0_1(v) \, dv <\infty.
$$
Then Proposition \ref{prop-mom4} in the appendix ensures that
$$
\sup_{\t \geq 0} \int_{\rr^3} \vert v \vert^4 \, g_1(\t, v) \, dv
<\infty,
$$
so that
$$
\sup_{\t \geq 0} \int_{\vert v \vert > R} \vert v \vert^2 \,
g_1(\t, v) \, dv
$$
tends to $0$ as $R$ goes to infinity. Prohorov's compactness
theorem and Proposition \ref{w2properties} imply the existence of
a sequence $\t_k\to\infty$ as $k\to\infty$ and a probability
measure $\mu^0$ on $\rr^3$ with zero mean velocity and unit
kinetic energy such that $W_2(g_1(\t_k), \mu^0)\to 0$ as
$k\to\infty$. We want to prove that $\mu^0 = g_{\infty}$.

\medskip

Without loss of generality, we can assume that the diverging time
sequence satisfies $\t_k + 1 \leq \t_{k+1}$ for all $k$. Now,
since $g_{\infty}$ is a stationary solution, it follows from the
first part of the theorem that
\begin{equation}\label{leq3}
W_2(g_1({\t_{k+1}}), g_{\infty}) \leq W_2(g_1(\t_{k} +1),
g_{\infty}) \leq W_2(g_1({\t_{k}}), g_{\infty}).
\end{equation}
On one hand, both $W_2(g_1({\t_{k}}), g_{\infty})$ and
$W_2(g_1({\t_{k+1}}), g_{\infty})$ tend to $W_2(\mu^0,
g_{\infty})$ as $k$ goes to infinity by triangular inequality.
Then, if $\mu(\t)$ denotes the solution to \eqref{pde-g1} with
initial datum $\mu^0$, the first point again ensures that
$$
W_2(g_1(\t_{k} +1), \mu (1)) \leq W_2(g_1(\t_k), \mu^0)
$$
which tends to $0$. Hence $W_2(g_1(\t_{k} +1), g_{\infty})$ tends
to $W_2(\mu(1), g_{\infty})$ by triangular inequality, and finally
$$
W_2(\mu(1), g_{\infty}) = W_2(\mu^0, g_{\infty})
$$
by passing to the limit in $k$ in \eqref{leq3}. By the non-increasing
character of $W_2$ along the flow, we deduce that
$$
W_2(\mu(1), g_{\infty}) = W_2(\mu(\t), g_{\infty})=W_2(\mu^0,
g_{\infty})
$$
for all $\t\in[0,1]$.

\medskip

Consequently $\mu(\t)$ and $g_{\infty}$ are two solutions to
\eqref{pde-g1} with zero mean velocity and unit temperature, whose
$W_2$ distance is constant on the time interval $[0, 1]$. This is
possible only if equality holds at each step in the proof of
Theorem \ref{thm-BE1tau} in the original space variables; in
particular
$$
W_2(Q^+(\mu(\t), \mu(\t)), Q^+(g_{\infty}, g_{\infty})) =
\sqrt{\frac{3+e^2}{4}}W_2(\mu(\t), g_{\infty})
$$
for all $\t$, and especially for $\t = 0$. But $\mu^0$ and $g_{\infty}$ have same mean
velocity and temperature, and, according to
\cite[Theorem 5.3]{Bobylev-Cercignani2}, $g_{\infty}$ is
absolutely continuous with respect to Lebesgue measure, with
positive density. Hence Proposition \ref{eg-Q+} ensures that
$\mu^0 = g_{\infty}$.

\medskip

In particular $W_2(g_1(\t_k), g_{\infty}) \to 0$ as $k\to \infty$,
and then $W_2(g_1(\t), g_{\infty})\to 0$ as $\t\to \infty$ since it
is a non increasing function.

\bigskip

{\bf Step 2.-} Let us now remove the hypothesis on the boundedness
of the initial fourth order moment. Let $(g^{0n})_n$ be a sequence in
$\Ptwo$ with zero mean velo\-city, unit kinetic energy, finite
fourth order moment and converging to $g_1^0$ in the weak sense of
probability measures; in particular it converges to $g_1^0$ in the
$W_2$ distance sense by Proposition \ref{w2properties}. Such a
$g^{0n}$ can be obtained by successive truncation of $g_1^0$ to a
ball of radius $n$ in $\rr^3$, translation to keep the mean
property, and dilation centered at $0$ to keep the kinetic energy
equal to $1$.

\smallskip

Then, if $g^n(\t)$ is the solution to \eqref{pde-g1} with initial
datum $g^{0n}$, the triangular inequality for $W_2$ and
\eqref{contr1g} ensure that
\begin{align*}
W_2(g_1(\t), g_{\infty}) &\leq W_2(g_1(\t), g^{n}(\t)) +
W_2(g^n(\t), g_{\infty}) \\
& \leq W_2(g_1^0, g^{0n}) + W_2(g^n(\t), g_{\infty}).
\end{align*}
Given $\eps >0$, the first term in the right hand side is bounded
by $\eps$ for some $n$ large enough, and for this now fixed $n$,
the second term is bounded by $\eps$ for all $\t$ larger than some
constant by the first step. This ensures that $W_2(g_1(\t),
g_{\infty})$ tends to $0$ as $\t$ goes to infinity.
\endproof

\bigskip

\subsection{The diffusive case}\label{sec-diff}

We now turn to the diffusive version \eqref{BE2} of \eqref{BE1}.
Again by the change of time
\[
\t = \frac{B}{E} \int_0^t \sqrt{\theta (f(w))} \, dw
\]
with $\displaystyle E = \frac{8}{1-e^2}$
we are brought to studying the equation
\begin{equation}\label{BE2tau}
\frac{\p f}{\p \t}= E \, Q(f,f) + \Theta^2(f(\tau)) \, \Delta_v f
 \end{equation}
where
$$
\qquad \Theta^2(f(\t)) = \frac{E \, A}{B} \, \big[\theta (f(\t)) \big]^{p-1/2}.
$$

\bigskip

As in the nonviscous case of \eqref{BE1tau} we shall prove

\begin{theorem}\label{thm-BE2tau}
If $f_1$ and $f_2$ are two solutions to \eqref{BE2tau} for the
respective initial data $f_1^0$ and $f_2^0$ in $\Ptwo$ with same kinetic energy, then
\begin{align}
W_2^2(f_1(\t), f_2(\t)) \leq {\rm e}^{-2 \tau}\, W_2^2(f_1^0,
f_2^0) + (1 - {\rm e}^{-2 \tau})\,\big\vert \!<\!f_1^0\!> -
<\!f_2^0\!> \!\big\vert^2 \label{cont-BE2tau}
\end{align}
for all $\t \geq 0$.
\end{theorem}

\proof We again start by giving a Duhamel's representation of the
solutions. To this aim we write \eqref{BE2tau} as
\[
\frac{\p f}{\p \t}= E \, F - E \, f+ \Theta^2(f(\tau)) \, \Delta f
\]
where $F = Q^+(f, f)$,  that is,
\[
\frac{\p \hat{f}}{\p \t} + \big(E + \vert k \vert^2 \, \Theta^2 (f(\t)) \big)
\, \hat{f} = E \, \hat{F}.
\]
Here, we are using the convention
\[
\hat{\mu}(k) = \int_{\rr^3} {\rm e}^{-i \, k \cdot x} \, d \mu( x)
\]
for the Fourier transform of the measure $\mu$ on $\rr^3$. Hence
the solutions satisfy
\[
\hat f (\t, k) = {\rm e}^{-E \, \t} \, \hat f^0(k) \, {\rm
e}^{-\Sigma(f, \t) \, \vert k \vert^2} + E \int_0^{\t} {\rm e}^{-E
(\t-s)} \, \hat F (s, k) \, {\rm e}^{-(\Sigma(f, \t) - \Sigma(f,
s)) \vert k \vert^2} \, ds
\]
where $\ds \Sigma (f, \t) = \int_0^{\t} \Theta^2 (f(s)) \, ds $, and
thus
\begin{align*}
f(\t, v) &= {\rm e}^{-E \, \t} \,  (f^0 * \Gamma_{2 \Sigma(f,
\t)}) (v) + E \! \int_0^{\t} \! \! {\rm e}^{-E (\t-s)} \, (F(s) *
\Gamma_{2 (\Sigma(f, \t) - \Sigma(f, s))})(v) \, ds\\
&:= {\rm e}^{-E \, \t} \,  \tilde{f}(\t, v) + E \int_0^{\t} {\rm
e}^{-E (\t-s)} \, \tilde{F}(\t,s, v) \, ds.
\end{align*}
Here
$$
\Gamma_{\alpha}(v) = \frac{1}{(2 \pi \alpha)^{3/2}} \, {\rm e}^{-
\vert v \vert^2 / 2 \alpha}
$$
is the centered Maxwellian with temperature $\alpha/3 >0$.
Moreover $f_1$ and $f_2$ have same temperature at all times, so
that $\Sigma(f_1, \t) = \Sigma(f_2, \t)$. Then the convexity of
the squared Wasserstein distance and its non-increasing character
by convolution with a given measure, see
Proposition \ref{w2properties}, imply that
\begin{align*}
W_2^2(f_1(\t), f_2(\t)) \!&\leq \!{\rm e}^{-E \t}
W_2^2(\tilde{f}_1(\t),
\tilde{f}_2(\t)) \!+\! E \!\!\int_0^{\t} \!\!\!{\rm e}^{-E(\t-s)} W_2^2(\tilde{F}_1(\t,s),\tilde{F}_2(\t,s)) ds \\
& \leq {\rm e}^{-E \t} \, W_2^2(f_1^0, f_2^0) + E \int_0^{\t} {\rm
e}^{-E(\t-s)} \, W_2^2 (F_1(s), F_2(s)) \, ds.
\end{align*}
In other words the squared distance $W_2^2(f_1(\t), f_2(\t))$
satisfies the same bound as in the nonviscous case of Theorem
\ref{thm-BE1tau}, and we can conclude analogously.
\endproof

\medskip

\begin{remark}
\
\begin{enumerate}
\item As pointed out to us by C. Villani the result can also be
obtained by a splitting argument between the collision term and
the diffusion term.

\item As proven in {\rm \cite{BisiCT}}, the temperature
$\theta(f(t))$ of the solution $f$ in the original time variable $t$
converges towards
\[
\theta_{\infty} = \Big( \frac{8 \, A}{B (1 - e^2)}
\Big)^{\frac{2}{3 - 2p}}
\]
as $t$ goes to infinity, and satisfies $\theta(f(t)) \, \geq \, \min(
\theta (f(0)), \theta_{\infty})$. In particular
\[
\t = \frac{B}{E} \int_0^t \sqrt{\theta (f(s))} \, ds \geq
\frac{C_1}{E} \, t
\]
if $C_1 = B \, \min( \theta (f(0)), \theta_{\infty})^{1/2}$. Writing
\eqref{cont-BE2tau} in the original variable $t$ for initial data
with equal mean velocity and temperature, we recover the contraction property
$$
W_2(f_1(t),f_2(t)) \leq W_2 (f_1^0,f_2^0)\, \rm{e}^{-(1-\gamma)C_1
t} ,
$$
that coincides with {\rm (3.1)} in {\rm \cite{BisiCT}} for the
Fourier-based $d_2$ distance exactly with the same rate. For $p=1$
one can exactly compute $\t$ and also recover {\rm (3.2)} in {\rm
\cite{BisiCT}} but for the distance $W_2$.

\item The existence of unique diffusive equilibria for each given
value of the initial mean velocity can be obtained from
this contraction property of the $W_2$ distance analogously to the
arguments done in {\rm \cite{BisiCT}} with the Fourier-based
distance $d_2$.
\end{enumerate}
\end{remark}


\medskip

\section{General cross section}\label{sec-gene-b}

In this section, we consider the more general case of a variable
collision cross section when the gain term $Q^+$ is defined by
$$
 (\varphi,Q^+(f,f)) = \frac{1}{4 \pi} \int_{\Realr^3} \int_{\Realr^3}
\int_{S^2} f(v) \, f(w) \, \varphi (v') \, b \Big(\frac{v-w}{\vert
 v-w\vert} \cdot \sigma\Bigl) \,d\sigma\,dv\,dw
$$
where again the post-collisional velocity $v'$ is given by
$$
 \dis v'= \frac{1}{2} (v+w)+ \frac{1-e}{4}(v-w)+\frac{1+e}{4} |v-w|\, \sigma
$$
and the cross section $b$ satisfies the normalized cut-off assumption
\begin{equation}\label{cutoff}
\int_{S^2} b(k \cdot \sigma) \, d\sigma = \int_0^{2 \pi}
\int_0^{\pi} b(\cos \theta) \, \sin \theta \, d\theta
\, d\phi = 1
\end{equation}
for any $k$ in $S^2$. Then we shall prove the
following extension of Proposition \ref{prop-qplus} for non
constant cross sections $b$:

\begin{theorem}\label{thm-b}
If $f$ and $g$ in $\Ptwo$ have equal mean velocity, then
\[
W_2^2(Q^+(f, f) , Q^+(g, g) ) \!\leq \!\left(\!\frac{3+ e^2}{4} \!+\!
\frac{1-e^2}{2} \pi\! \int_0^{\pi}\! \!b(\cos \theta) \cos \theta \,
\sin \theta \, d\theta \right) W_2^2(f, g).
\]
\end{theorem}

\medskip

Before going onto the proof, we draw the main consequence. Let
$f_1=f_1(\t,v)$ and $f_2=f_2(\t,v)$ be two solutions to the
Boltzmann equation
$$
\frac{\partial f}{\partial \t} = Q(f, f) = Q^+(f, f) - f
$$
with respective initial data $f_1^0$ and $f_2^0$ in $\Ptwo$, where
$Q^+$ is defined as above. Then, as in Section \ref{sec-nodiff},
Duhamel's representation formula
$$
f(\t) = {\rm e}^{-\t} \, f^0 +  \int_0^{\t} {\rm e}^{-(\t-s)} \,
Q^+(f(s), f(s)) \, ds
$$
of the solutions and the convexity of $W_2^2$ ensure the contraction property
\begin{equation}\label{cont-edp-b}
W_2 (f_1(\t), f_2(\t)) \leq {\rm e}^{-(1- \gamma_b) \t /2} \,  W_2
(f_1^0, f_2^0)
\end{equation}
for all $\t$, where
$$
\gamma_b = \frac{3+ e^2}{4} + \frac{1-e^2}{2} \,\pi \int_0^{\pi}
b(\cos \theta) \, \cos \theta \, \sin \theta \, d\theta
$$
is bounded by $1$ by \eqref{cutoff}.

\medskip

In the elastic case when $e=1$, $\gamma_b =1$, one recovers
Tanaka's non-strict contraction  result \cite{Tanaka} for the
solutions to the homogeneous elastic Boltzmann equation for
Maxwellian molecules, at least under the cut-off assumption, but
with a somehow simpler argument than those given in \cite{Tanaka}
and \cite{Vil03}.

\

\proof By definition
\begin{eqnarray*}
(\varphi,Q^+(f,f))
& =&
2 \pi \int_0^{\pi} \!\int_{\Realr^3} \int_{\Realr^3} \Big\{ \int_0^{2 \pi} \varphi
(v') \, \frac{d\phi}{2\pi} \Big\}  f(v) \, f(w) \, dv \, dw  \,b(\cos \theta)
\, \sin \theta \, d\theta \\
& = &
2 \pi \int_0^{\pi} \ee \, \big[ (\varphi, {\mathcal U}_{V, W, \theta}) \big]\,  b(\cos \theta)
\, \sin \theta \, d\theta
\end{eqnarray*}
where $V$ and $W$ are independent random variables distributed
according to $f$ and, given $v, w$ in $\rr^3$, ${\mathcal U}_{v, w, \theta}$
is the uniform probability measure on the circle $C_{v, w,
\theta}$ with center
\[
c_{v, w, \theta} = \frac{1}{2} (v+w)+ \Big(\frac{1-e}{4} +
\frac{1+e}{4} \cos\theta \Big) (v-w),
\]
radius
\[
r_{v, w, \theta} = \frac{1+e}{4} |v-w| \sin\theta
\]
and axis
\[
k = \frac{v-w}{\vert v-w \vert} \cdot
\]

\begin{figure}[ht]
\centering
\includegraphics[scale=0.7]{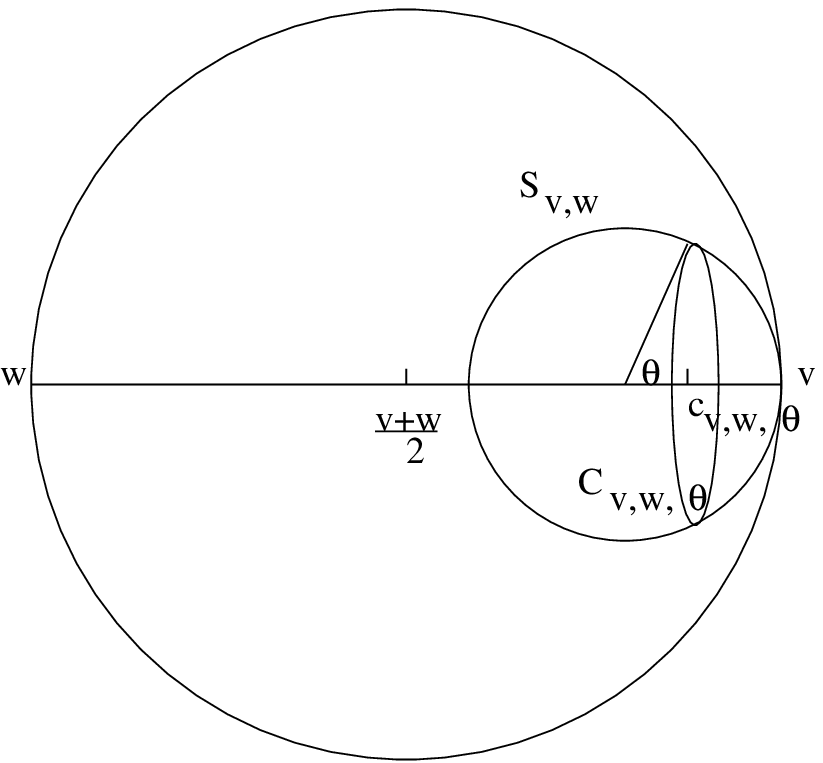}
\caption{} \label{sphere-circle}
\end{figure}

Let also $g$ be a Borel probability measure on $\rr^3$ and $X, Y$
be independent random variables with law $g$. Then, by the
normalization assumption \eqref{cutoff}, the convexity of the
squared Wasserstein distance with respect to both arguments
ensures that
\begin{equation}\label{conv-b}
W_2^2(Q^+(f, f) , Q^+(g, g) ) \leq
2 \pi \!\int_0^{\pi} \! \ee \,
\big[ W_2^2({\mathcal U}_{V, W, \theta}, {\mathcal U}_{X, Y, \theta}) \big]\,  b(\cos \theta)
\, \sin \theta \, d\theta.
\end{equation}

We now let $v, w, x, y$ and $\theta$ be fixed in $\rr^3$ and $[0,
\pi]$ respectively, and give an upper bound to $W_2^2({\mathcal U}_{v, w,
\theta}, {\mathcal U}_{x, y, \theta})$. This consists in estimating the
transport cost of a circle in $\rr^3$ onto another one, for which
we have the following general bound:

\begin{lemma}{\rm \cite{Vil03}}\label{w2circle}
The squared Wasserstein distance between the uniform distributions on  the circles with centers $c$
and $c'$, radii $r$ and $r'$ and axes $k$ and $k'$ is bounded by
$$
\vert c - c' \vert^2 + r^2 + r'^2 - r r' (1 + \vert k \cdot k' \vert).
$$
\end{lemma}

\noindent Hence, using the notations $a=v-x$, $b=w-y$,
$\tilde{a}=v-w$ and $\tilde{b}=x-y$ in our case we get
\begin{align}\label{boundUtheta}
W_2^2({\mathcal U}_{v, w,\theta}, {\mathcal U}_{x, y, \theta}) \leq & \;
\Big\vert \Big(
\frac{3-e}{4} + \frac{1+e}{4} \cos \theta \Big) a + \frac{1+e}{4}
(1 - \cos \theta) \, b \Big\vert^2 \nonumber
\\
 & \, + \Big(\frac{1+e}{4}\Big)^2 \sin^2 \theta \Big[ \vert \tilde{a}
\vert^2 +  \vert \tilde{b} \vert^2 -  \vert \tilde{a} \vert \vert
\tilde{b} \vert \Big(1 + \Big( \frac{\tilde{a}}{\vert \tilde{a}
\vert} \cdot \frac{\tilde{b}}{\vert \tilde{b} \vert} \Big) \Big]
\nonumber
\\
\leq & \; \Big[\Big( \frac{3-e}{4} + \frac{1+e}{4} \cos \theta
\Big)^2 +
  \Big(\frac{1+e}{4}\Big)^2 \sin^2 \theta \Big] \vert a \vert^2
  \nonumber
\\
 & \, + \, 2 \, \Big(\!\frac{1+e}{4}\!\Big)^2 \Big[ \Big( \frac{3-e}{1+e} + \cos
\theta \Big) (1 - \cos \theta) -2 \sin^2 \theta \Big] a \cdot b \nonumber
\\
 & \, + \Big(\frac{1+e}{4}\Big)^2 \big[ (1 - \cos \theta)^2 + \sin^2
\theta \big] \vert b \vert^2
\end{align}
where we have used the bound
\begin{multline*}
\vert \tilde{a}\vert^2 +  \vert \tilde{b}\vert^2 -  \vert
\tilde{a}\vert \vert \tilde{b}\vert - \tilde{a} \cdot \tilde{b}
\leq \vert \tilde{a}\vert^2 +  \vert \tilde{b}\vert^2 - 2\,
\tilde{a} \cdot \tilde{b} = \vert\tilde{a} - \tilde{b} \vert^2 =
\vert a - b \vert^2.
\end{multline*}
Assume now that $(V, X)$ and $(W, Y)$ are two independent couples of random
variables, optimal in the sense that
$$
W_2^2(f,g) = \ee \, \big[ \vert V-X  \vert^2 \big] = \ee \, \big[ \vert W-Y  \vert^2 \big].
$$
Note that
$$
\ee \, \big[ (V-X) \cdot (W-Y) \big] = \ee \, \big[ (V-X) \big] \cdot \ee \, \big[ (W-Y) \big]  = 0
$$
since $(V, X)$ and $(W, Y)$ are independent and since $f$ and $g$ have
same mean velocity. Then from \eqref{boundUtheta}:
$$
\ee \, \big[ W_2^2({\mathcal U}_{V, W, \theta}, {\mathcal U}_{X, Y, \theta}) \big]  \leq
\gamma(\theta) \, W_2^2(f, g)
$$
where
\begin{eqnarray*}
\gamma (\theta) &=& \Big( \frac{3-e}{4} + \frac{1+e}{4} \cos
\theta \Big)^2 +
    \Big(\frac{1+e}{4}\Big)^2 \big[ (1 - \cos \theta)^2 + 2\sin^2
      \theta \big] \\
&=&
\frac{3+e^2}{4} + \frac{1-e^2}{4} \cos\theta.
\end{eqnarray*}
One concludes the argument after averaging over $\theta$ as in
\eqref{conv-b} and taking \eqref{cutoff} into account.
\endproof


\section{Inelastic Kac Model}\label{sec-kac}

In this last section we consider a simple one-dimensional model
introduced in \cite{PT03} which can be seen as a dissipative
version of the Kac caricature of a Maxwellian gas
\cite{Kac,McKean}. Let us remark that the definition and
properties of the Euclidean Wasserstein distance $W_2$ discussed
above generalizes equally well to any dimension. Tanaka himself
\cite{Tanaka73} showed that the Euclidean Wasserstein distance is a
non strict contraction for the elastic classical Kac model. In the
inelastic Kac model, the evolution of the density function $f$ is
governed by the equation
\begin{equation}\label{eq-kac}
\frac{\partial f}{\partial t} = Q(f, f)
\end{equation}
in which the collision term $Q(f,f)$ is defined by
$$
(\varphi, Q(f,f)) = \int_{\Realr} \int_{\Realr}
\int_0^{2 \pi} f(v) f(w) \Big[ \varphi (v')- \varphi (v) \Big]
 \, \frac{d\theta}{2\pi}\,dv\,dw
$$
for any test function~$\varphi$, where
\[
 \dis v'= v \, \cos \theta \,  \vert \cos  \theta \vert^p -  w \, \sin
 \theta \,  \vert \sin \theta \vert^p
\]
is the postcollisional velocity and $p > 0$ measures the
inelasticity. Equation \eqref{eq-kac} preserves mass but makes the
momentum and kinetic energy decrease to $0$ at an exponential
rate, $\theta(f(t))={\rm e}^{-2\beta t}\theta (f^0) + ({\rm
e}^{-2\beta t} -{\rm e}^{-2 t}) <\!\!f^0\!\!>$ with $\beta > 0$
given below. In particular, solutions to \eqref{eq-kac} tend to
the Dirac mass at $0$.

As in the inelastic Maxwell model discussed above, we start by
deriving a contraction property for the gain operator $Q^+$ defined by
$$
(\varphi, Q^+(f,f)) = \int_{\Realr} \int_{\Realr}
\int_0^{2 \pi} f(v) f(w) \, \varphi (v') \,\frac{d\theta}{2\pi}\,dv\,dw .
$$

\begin{proposition}\label{propkac}
If $f$ and $g$ belong to ${\mathcal P}_2(\rr)$, then
$$
W_2^2(Q^+(f,f), Q^+(g,g)) \leq \left[\int_0^{2 \pi} \big( \vert \cos \theta
\vert^{2(p+1)} + \vert \sin \theta \vert^{2(p+1)}\big)  \, \frac{d\theta}{2\pi}
\right]\, \, W_2^2(f,g).
$$
\end{proposition}

\medskip

In terms of solutions $f(t)$ and $g(t)$ to the modified Kac equation
\eqref{eq-kac} with finite initial energy only, the above proposition yields, as in previous sections,
the bound
$$
W_2(f(t), g(t)) \leq {\rm e}^{- \beta t} \, W_2(f^0, g^0)
$$
where
$$
2 \, \beta = 1 - \int_0^{2 \pi} \big(\vert \cos \theta
\vert^{2(p+1)} + \vert \sin \theta \vert^{2(p+1)}\big) \,
\frac{d\theta}{2\pi} \, > 0.
$$
This bound is optimal without further assumptions on the initial
data $f^0$ and $g^0$ since equality holds in the case when $<f^0>
= 0$ and $g^0=\delta_0$ analogously to previous cases.

\medskip

\proof
Given a vector $(v,w)$ in $\rr^2$, let $\mathcal C_{v,w}$ denote the curve
$$
\{(v'(\theta), w'(\theta)), \theta \in [0, 2\, \pi]\}
$$
where
\begin{equation}\label{defv'kac}
\left. \begin{array}{ccl}
v'(\theta) &=& v \, \cos \theta \, \vert \cos \theta \vert^p -  w \, \sin
\theta \, \vert \sin \theta \vert^p \\
w'(\theta) &=& v \, \sin \theta \, \vert \sin \theta \vert^p +  w \, \cos
\theta \, \vert \cos \theta \vert^p.
\end{array}
\right.
\end{equation}

Let also $\mathcal U_{v,w}$ be the uniform probability distribution
on  $\mathcal C_{v,w}$.

Given $V$ and $W$ two independent random variables distributed
according to $f$, we note that $Q^+(f,f)$
is the first marginal on $\rr$ of $ \ee \big[ \mathcal U_{V,W} \big]$,
but also its second marginal by symmetry.
Then, we have the following result, which is the analogous of Lemmas
\ref{w2unif} and \ref{w2circle} for this model:

\begin{lemma}\label{lemmakac}
Given two vectors $(v,w)$ and $(x,y)$ in $\rr^2$, the squared Wasserstein distance  between the
distributions
$\mathcal U_{v,w}$ and $\mathcal U_{x,y}$ is bounded by
$$
(1-2\beta)\, \big(\vert v-x \vert^2 + \vert w-y \vert^2 \big).
$$
\end{lemma}

\proof One can transport the curve
$\mathcal C_{v,w}$ onto $\mathcal C_{x,y}$ by the linear map
$$
(a,b) \mapsto T(a, b) = \frac{r'}{r} (a \, \cos \omega  - b \,
\sin \omega, a \, \sin \omega + b \, \cos \omega)
$$
where $r=\sqrt{v^2 + w^2}$, $r' =\sqrt{x^2 + y^2}$ and $\omega$ is
the angle between the vectors $(v,w)$ and $(x,y)$ in case they do
not vanish. We leave the reader  discuss the case when
either $(x,y)$ or $(v,w)$ are zero. Then, analogously to the proof of Lemma \ref{w2unif}, one can define a transport plan associated to the transport
map $T$ to get
$$
W_2^2(\mathcal U_{v,w}, \mathcal U_{x,y}) \leq \int_{\rr^2} \vert
T(a,b) - (a,b) \vert^2 \, d\mathcal U_{v,w} (a,b).
$$
Furthermore, for all $(a,b)$ in $\rr^2$,
\begin{eqnarray*}
 \vert T(a,b) - (a,b) \vert^2
&=& \Big\vert  \frac{r'}{r} \big( a \, \cos \omega  - b \, \sin
 \omega \big) -a \Big\vert^2 + \Big\vert  \frac{r'}{r} \big( a \, \sin \omega + b \, \cos
 \omega \big) -b \Big\vert^2 \\
&=& \Big( \Big(\frac{r'}{r} \Big)^2 - 2 \, \frac{r'}{r} \, \cos
\omega
 + 1 \Big) \big( a^2 + b^2 \big) \\
&=& \frac{\vert v-x \vert^2 + \vert w-y\vert^2}{v^2 + w^2} \, \big(
 a^2 + b^2 \big).
\end{eqnarray*}
Hence, we deduce
\begin{eqnarray*}
W_2^2(\mathcal U_{v,w},
\mathcal U_{x,y})
&\leq&
\frac{\vert v-x \vert^2
  + \vert w-y\vert^2}{v^2 + w^2} \int_{\rr^2}  (a^2 + b^2) \, d\mathcal U_{v,w}
  (a,b)\\
&=&
\frac{\vert v-x \vert^2  + \vert w-y\vert^2}{v^2 + w^2} \int_0^{2\pi} (v'(\theta)^2 +
w'(\theta)^2) \, \frac{d\theta}{2\pi} \cdot
\end{eqnarray*}
But
$$
v'(\theta)^2 + w'(\theta)^2 = \big( \vert \cos \theta
\vert^{2(p+1)} + \vert \sin \theta \vert^{2(p+1)} \big) (v^2 + w^2)
$$
by \eqref{defv'kac}, so that
$$
W_2^2(\mathcal U_{v,w}, \mathcal U_{x,y})
\leq
\left[ \int_0^{2 \pi} \big(\vert \cos \theta
\vert^{2(p+1)} + \vert \sin \theta \vert^{2(p+1)}\big) \,
\frac{d\theta}{2\pi}\right] \, \big( \vert v-x \vert^2  + \vert w-y\vert^2 \big)
$$
which is the bound given by the lemma.
\endproof

\

We now continue the {\it proof of Proposition \ref{propkac}}. First of
all, let $(V,X)$ and $(W,Y)$ be two independent couples of random
variables, with $V$ and $X$ distributed according to $f$, $W$ and
$Y$ according to $g$, optimal in the sense that
$$
W_2^2(f,g) = \ee \big[\vert V-W \vert^2 \big] = \ee \big[\vert X-Y
  \vert^2 \big].
$$
Then, by convexity of the squared Wasserstein distance again, it
follows from Lemma \ref{lemmakac} that
\begin{align}\label{w22kac}
W_2^2 \big( \ee \big[\mathcal U_{V,W} \big], \ee \big[ \mathcal
  U_{X,Y} \big] \big) & \leq \ee \big[ W_2^2 (\mathcal U_{V,W},
\mathcal U_{X,Y}) \big]\nonumber
\\
 &\leq (1-2\beta) \big(\ee \big[\vert V-W \vert^2 \big] + \ee \big[\vert X-Y
  \vert^2 \big] \big)\nonumber
\\
&= 2 \,(1-2\beta) \,  W_2^2(f,g).
\end{align}

Next, we remark that the measure $\mathcal U_{V,W}$ on $\rr^2$ has
first {\bf and} second marginals equal by symmetry of the curve
${\mathcal C}_{V,W}$ by a $\pi/2$ rotation. This implies that the
first and second marginals of $\ee \big[ \mathcal U_{V,W} \big]$
on $\rr^2$ are equal to $Q^+(f,f)$, and likewise for the
measure $\ee \big[ \mathcal U_{X,Y}\big]$ with marginals
$Q^+(g,g)$. We shall conclude the argument of Proposition \ref{propkac}
by using the following general result:

\begin{lemma}\label{lemmakac2}
If the Borel probability measures $\mu_j^i$ on $\rr$ are the successive
one-dimensional marginals of the measure $\mu^i$ on $\rr^N$, for $i=1, 2$ and $j=1,
\dots, N$, then
$$
\sum_{j=1}^N  W_2^2(\mu_j^1, \mu_j^2) \, \leq \, W_2^2(\mu^1, \mu^2).
$$
\end{lemma}

\proof Let $\pi$ be a measure on $\rr^N_v \times
\rr^N_w$ with marginals $\mu^1$ and $\mu^2$, optimal in the sense that
$$
W_2^2(\mu^1, \mu^2) =  \iint_{\rr^N \times \rr^N} \vert v -w \vert^2
\, d\pi(v, w).
$$
Then its marginal $\pi_j$ on $\rr_{v_j} \times \rr_{w_j}$ has itself
marginals $\mu_j^1$ and $\mu_j^2$, so
$$
W_2^2(\mu_j^1, \mu_j^2) \leq  \iint_{\rr \times \rr} \vert v_j -w_j \vert^2
\, d\pi_j(v_j, w_j).
$$
The lemma follows by noting that
$\displaystyle \sum_{j=1}^N  \vert v_j -w_j \vert^2 = \vert v -w \vert^2.$
\endproof

\

In our particular case, Lemma \ref{lemmakac2} ensures that
$$
2 \, W_2^2(Q^+(f,f), Q^+(g,g)) \leq
W_2^2 \big( \ee \big[\mathcal U_{V,W} \big], \ee \big[ \mathcal
  U_{X,Y} \big] \big)
$$
which concludes the proof of Proposition \ref{propkac} taking
\eqref{w22kac} into account.
\endproof


\section*{Appendix: Uniform in time Propagation of fourth order Moments}

In this appendix we derive a uniform propagation of fourth order
moments $\displaystyle \int_{\rr^3} \vert v \vert^4 \, g(\t, v) \, dv$
of solutions $g$ to
\begin{equation}\label{pde-g}
\frac{\partial g}{\partial \t} + \nabla \cdot (g \, v) = E \, Q(g, g)
\end{equation}
where the operator $Q(g, g)$ is defined as in \eqref{Qweak} for $0
< e < 1$ and $E = \displaystyle \frac{8}{1 - e^2} \cdot$

This result has been used in the proof of Theorem \ref{norate}.

\begin{proposition}\label{prop-mom4}
If $g^0$ is a Borel probability measure on $\rr^3$ such that
$$
\int_{\rr^3} \vert v \vert^4 \, g^0(v) \, dv<\infty,
$$
then the solution $g$ to \eqref{pde-g} with
initial datum $g^0$ verifies
$$
\sup_{\t \geq 0} \int_{\rr^3} \vert v \vert^4 \, g(\t, v) \,
dv<\infty.
$$
\end{proposition}

\proof Without loss of generality we can assume that $g^0$, and hence $g(\t)$
for all $ \t \geq 0$, has zero mean velocity. We let
$$
m_4(\t) = \int_{\rr^3} \vert v \vert^4 \, g(\t, v) \, dv
$$
denote the fourth order moment of $g(\t)$. Then, using the weak formulation of
the inelastic Boltzmann
equation, we have:
\begin{equation}\label{eq-mom4}
\frac{d m_4(\t)}{d \t} = \int_{\rr^3} \nabla (\vert v \vert^4)
\cdot v \, g(\t, v) \, dv + E \int_{\rr^3} \vert v \vert^4 \, Q(g(\t),
g(\t)) (v) \, dv.
\end{equation}
While the first term in the right hand side is simply $4 \,
m_4(\t)$, the second term is computed by

\begin{lemma}\label{lem-mom4}
There exist some constants $\mu_1$ and $\mu_2$, depending only on $e$, such that
\begin{align*}
\int_{\rr^3} \vert v \vert^4 \, Q(g, g) (v) \, dv = &- \lambda \, \int_{\rr^3} \vert v
\vert^4 \, g(v) \, dv   + \mu_1  \Big( \int_{\rr^3} \vert v
\vert^2 \,  g(v) \, dv \Big)^2 \\ &+ \mu_2 \iint_{\rr^3 \times \rr^3} (v \cdot
w)^2 \, g(v) \, g(w) \, dv \, dw
\end{align*}
for any probability measure $g$ on $\rr^3$ with finite moment of
order $4$ and zero mean velocity, where
$$
\lambda =  \frac{1}{3} ( 1 + 4 \, \eps - 7 \, \eps^2 + 4 \, \eps^3
- 2 \, \eps^4) \qquad \mbox{and} \qquad \eps = \frac{1-e}{2} \cdot
$$
\end{lemma}

\noindent With this lemma in hand, \eqref{eq-mom4} reads
\begin{equation}\label{new}
\frac{d m_4(\t)}{d \t} =  \Big( 4 - E \, \lambda \Big) m_4(\t) +
m(\t)
\end{equation}
where $m(\t)$ is a combination of second order moments, which are
bounded in time since the kinetic energy is preserved by equation
\eqref{pde-g}. Moreover one can check from the expression of $E$
and $\lambda$ in terms of $\eps = (1-e)/2$ that
$$
4 - E \, \lambda  = \frac{2}{3 \, \eps (1 - \eps)} \, [ -1 + 2 \,
\eps + \eps^2 - 4 \, \eps^3 + 2 \, \eps^4 ]
$$
which is negative for any $0 < \eps < 1/2$, that is, for any $0 < e <
1$. By Gronwall's lemma this ensures that $m_4(\t)$ is bounded
uniformly in time if initially finite, and concludes the argument to Proposition
\ref{prop-mom4}.
\endproof

\medskip

Let us remark that identity \eqref{new} is also useful to
understand that moments are not created by this equation in
contrast to the hard-spheres case \cite{MM,MMM}. In fact, if
initially moments are infinite, they will remain so. Thus, this is
another reason why homogeneous cooling states have only certain
number of moments bounded (see \cite{Bobylev-Cercignani2}).

\medskip

We now turn to the {\it proof of Lemma \ref{lem-mom4}}, whose
result is given in \cite[Section~4]{Bobylev-Cercignani} only in
the radial isotropic case, i.e., whenever $g(v)$ depends only on
$\vert v \vert$. By symmetry we start by writing
$$
\int_{\rr^3} \! \! \vert v \vert^4 \, Q(g, g) (v) \, dv =
\frac{1}{4 \pi} \int_{\Realr^3} \int_{\Realr^3}
\int_{S^2} \! \! g(v) g(w)  \frac{1}{2} [ \vert v' \vert^4 + \vert w' \vert^4 - \vert v \vert^4 - \vert w \vert^4 ]
\, d\sigma\,dv\,dw
$$
where
\begin{eqnarray*}
v' &=&  \frac{1}{2} (v+w) + \frac{1-e}{4}(v-w)+\frac{1+e}{4} |v-w|\sigma\\
w'&=& \frac{1}{2} (v+w) - \frac{1-e}{4}(v-w) - \frac{1+e}{4} |v-w|\sigma.
\end{eqnarray*}
Then we introduce the notation
$$
u = \frac{v+w}{2} \raise2pt \hbox{,} \quad U = \frac{v-w}{2} \raise2pt \hbox{,}
 \quad \eps = \frac{1-e}{2} \raise2pt \hbox{,} \quad \eps' = 1 - \eps = \frac{1+e}{2}
$$
in which
$$
v' = u + \eps \, U + \eps'  \, \vert U \vert \, \sigma, \quad v' =
u - \eps \, U - \eps'  \, \vert U \vert \, \sigma, \quad v = u +
U, \quad w = u - U.
$$
Then
\begin{eqnarray*}
\vert v' \vert^2 &=&  \vert u \vert^2 + (\eps^2  + \eps'^2) \vert U \vert^2 + 2\,  \eps \, \eps' \, \vert U \vert \, ( U \cdot \sigma) + 2 \, \eps \, ( u \cdot U) + 2 \, \eps' \, \vert U \vert \, (u \cdot \sigma) \\
\vert w' \vert^2&=& \vert u \vert^2 + (\eps^2  + \eps'^2) \vert U \vert^2 + 2 \, \eps \, \eps' \, \vert U \vert ( U \cdot \sigma) - 2 \, \eps \, (u \cdot U) - 2 \, \eps' \, \vert U \vert \, (u \cdot \sigma) \\
\vert v \vert^2 &=& \vert u \vert^2 + \vert U \vert^2 + 2 \, (u \cdot U) \\
\vert w \vert^2 &=& \vert u \vert^2 + \vert U \vert^2 - 2 \, (u \cdot U)
\end{eqnarray*}
and eventually
\begin{align*}
\frac{1}{2} [ \vert v' \vert^4 + \vert w' \vert^4 &- \vert v \vert^4 - \vert w \vert^4 ] \\
  = &\big[(\eps^2 + \eps'^2)^2 -1 \big] \vert U \vert^4 + 2(\eps^2 + \eps'^2 -1) \vert u \vert^2 \vert U \vert^2 + 4 (\eps^2 - 1) (u \cdot U)^2 \\
 & + 4 \, \eps^2 \, \eps'^2 \, \vert U \vert^2 \, (U \cdot \sigma)^2 + 4 \, \eps'^2 \, \vert U \vert^2 \, (u \cdot \sigma)^2 \\
 &+ 4 \, \eps \, \eps' \, \vert U \vert \, \big[ \vert u \vert^2 + (\eps^2 + \eps'^2) \vert U \vert^2 \big] \, (U \cdot \sigma)
  + 8 \, \eps \, \eps' \, \vert U \vert \, (u \cdot U) \, (u \cdot \sigma).
\end{align*}
Integrating with respect to $\sigma$ in $S^2$ and taking the
identities
$$
 \int_{S^2} 1 \, \frac{d \sigma}{4 \pi} = 1, \qquad  \int_{S^2} (k \cdot \sigma) \, \frac{d \sigma}{4 \pi} = 0, \qquad   \int_{S^2} (k \cdot \sigma)^2 \, \frac{d \sigma}{4 \pi} = \frac{\vert k \vert^2}{3}
$$
into account, we obtain
$$
 \int_{S^2}  \frac{1}{2} [ \vert v' \vert^4 + \vert w' \vert^4 - \vert v \vert^4 - \vert w \vert^4 ] \, \frac{d \sigma}{4 \pi}
 =  \alpha \, \vert U \vert^4 + \beta \,  \vert u \vert^2 \vert U \vert^2 +
  \gamma \, (u \cdot U)^2
$$
where
$$
\alpha = (\eps^2 + \eps'^2)^2 -1 + \frac{4}{3} \eps^2 \eps'^2, \quad  \beta = 2 \big[ \eps^2 + \eps'^2 -1 + \frac{2}{3} \eps'^2 \big], \quad \gamma = 4 (\eps^2 - 1).
$$

Then, by definition of $u$ and $U$ in terms of $v$ and $w$, the identities
$$
\iint_{\rr^3 \times \rr^3} \vert U \vert^4 \, g(v) \, g(w) \, dv \, dw =
 \frac{1}{8} [ m_4 + m_2^2 + 2 \, \overline{m_2^2}],
$$
$$
\iint_{\rr^3 \times \rr^3} \vert u \vert^2 \, \vert U \vert^2 \, g(v) \, g(w) \, dv \, dw =
\frac{1}{8} [ m_4 +  m_2^2  - 2 \,  \overline{m_2^2}]
$$
and
$$
\iint_{\rr^3 \times \rr^3} ( u \cdot U )^2 \, g(v) \, g(w) \, dv \, dw = \frac{1}{8} [m_4 - m_2^2]
$$
hold with
$$
m_4 = \int_{\rr^3} \vert v \vert^4 \, g(v) \, dv, \; \;
 m_2 = \int_{\rr^3} \vert v \vert^2 \, g(v) \, dv, \; \;
\overline{m_2^2} = \iint_{\rr^3 \times \rr^3} ( v \cdot w )^2 \, g(v) \, g(w) \, dv \, dw
$$
since $g$ has zero mean velocity. Collecting all terms, we obtain
$$
\int_{\rr^3} \vert v \vert^4 \, Q(g, g) (v) \, dv =
- \lambda \, m_4 + \mu_1 \, m_2^2 + \mu_2 \, \overline{m_2^2}
$$
where
$$
\lambda = - \frac{1}{8} (\alpha + \beta + \gamma) = \frac{1}{3} ( 1
+ 4 \, \eps - 7 \, \eps^2 + 4 \, \eps^3 - 2 \, \eps^4),
$$
$$
\mu_1 = \frac{1}{8} (\alpha + \beta - \gamma) \quad {\textrm{and}}
\quad \mu_2 = \frac{1}{4} (\alpha - \beta)
$$
depend only on $\eps$, that is, only on $e$. This concludes the proof of Lemma \ref{lem-mom4}.
\endproof


\

\noindent {\bf Acknowledgements.-} The authors are grateful to
Laurent Desvillettes, Giusep\-pe Toscani and C\'edric Villani for
stimulating discussions and fruitful comments. JAC acknowledges
the support from DGI-MEC (Spain) project MTM2005-08024.

\end{document}